\documentclass[12pt]{article}
\usepackage[scaled]{helvet} % ss
\usepackage{courier} % tt
\usepackage{euler} % math
\normalfont
\usepackage[T1]{fontenc}
\usepackage{etex} 
\usepackage{amsmath,geometry,amssymb,mathtools,latexsym,amsthm,enumerate,leftidx,tikz,url}
\usepackage{bussproofs}
\usepackage{listings,epigraph}
\usepackage{relsize}
\usepackage[all,cmtip]{xy}
\xyoption{2cell}
\xyoption{curve}
\UseTwocells
\input{diagxy}
\CompileMatrices      
\usepackage{tikz}
\usepackage{pdfpages}

\newcommand{\toarrow}{\ensuremath{\rightarrow}}

\newcommand{\Id}{\mathsf{Id}}
\newcommand{\myemph}[1]{\textbf{#1}}    % Produces boldface text
\newcommand{\imp}{\ensuremath{\Rightarrow}} 
\newcommand{\arr}{\ensuremath{\rightarrow}}

\newcommand{\U}{\mathcal{U}}      
\newcommand{\Sn}{\mathbb{S}}       
\newcommand{\Z}{\mathbb{Z}}       
    
\newcommand{\lloop}{\mathsf{loop}}    
\newcommand{\base}{\mathsf{base}}   
\newcommand{\cov}{\mathsf{cov}}   
\newcommand{\suc}{\mathsf{succ}}   
\newcommand{\ua}{\mathsf{ua}}    
       
\newcommand{\id}[1]{\mathsf{Id}_{#1}} 
\newcommand{\judge}[3][]{#2\;\vdash_{#1}\;#3}

\newtheorem{theorem}{Theorem}
\newtheorem*{theorem*}{Theorem}

\theoremstyle{remark}
 
\newtheorem*{remarks*}{Remarks}

\theoremstyle{definition}

\newtheorem*{definition*}{Definition}

\begin{document}
%%%%%%%%%%%%%%%%%%%%%%%%%%%%%%%%%%%%%%%%%%%%%%%%%%%%%%%%%%%%

\title{
A proposition is the (homotopy) type of its proofs}
\author{
Steve Awodey\thanks{
Thanks to Ulrik Buchholtz and Michael Shulman for comments on an earlier draft.  This research was partially supported by the U.S. Air Force Office of Scientific Research through MURI grant FA9550-15-1-0053. Any opinions, findings and conclusions or recommendations expressed in this material are those of the authors and do not necessarily reflect the views of the AFOSR.}\\
}
\date{
\today
}

\maketitle

\epigraph{
There are, at first blush, two kinds of construction involved: constructions of proofs of some proposition and constructions of objects of some type.  But I will argue that, from the point of view of foundations of mathematics, there is no difference between the two notions. A proposition may be regarded as a type of object, namely, the type of its proofs.  Conversely, a type $A$ may be regarded as a proposition, namely, the proposition whose proofs are the objects of type $A$.  So a proposition $A$ is true just in  case there is an object of type $A$.
}{W.W.~Tait \cite{Tait}}

%
%%%%%%%%%%%%%%%%%%%%%%%%%%%%%%%%%%%%%%%%%%%%%%%%%%%%%%%%%
\section*{Overview}
%%%%%%%%%%%%%%%%%%%%%%%%%%%%%%%%%%%%%%%%%%%%%%%%%%%%%%%%%

\emph{Homotopy type theory} is a new field devoted to a recently discovered connection between Logic and Topology --- more specifically, between constructive type theory, which was originally invented as a constructive foundation for mathematics and now has many applications in the theory of programming languages and formal proof verification, and homotopy theory, a branch of algebraic topology devoted to the study of continuous deformations of geometric spaces and mappings.
 The basis of homotopy type theory is an interpretation of the system of intensional type theory into abstract homotopy theory.  As a result of this interpretation, one can construct new kinds of models of constructive logic and study that system semantically, e.g.\ proving consistency and  independence results.  Conversely, constructive type theory can also be used as a formal calculus to reason about abstract homotopy.  This is particularly interesting in light of the fact that the type theory used underlies several computational proof assistants, such as Coq and Agda; this allows one to use those systems to reason formally  about homotopy theory and fully verify the correctness of definitions and proofs using these computer proof systems.  Potentially, this could provide a useful tool for mathematicians working in fields like homotopy theory and higher category theory.  Finally, new logical principles and constructions based on homotopical and higher categorical intuitions can be added to the system, providing a way to formalize many classical spaces and sophisticated mathematical constructions.  Examples include the so-called \emph{higher inductive types} and the \emph{univalence axiom} of Voevodsky.

More broadly, \emph{univalent foundations} is an ambitious new program for foundations of mathematics, proposed by Voevodsky, which is based roughly on homotopy type theory and intended to capture a very broad range of mathematics (I am reluctant to use the phrase ``All of Mathematics'', but  there is nothing in particular that could not, in principle, be done).   The new univalence axiom, which roughly speaking implies that isomorphic structures can be identified, and the general point of view that it promotes sharpen the expressiveness of the system and make it more powerful, so that new concepts can be isolated and new constructions can be carried out, and others that were previously ill-behaved (such as quotients) can be better controlled. The system is not only more expressive and powerful than previous type- and set-theoretic systems of foundations; it also has two further, distinct novelties: it is still amenable to computer formalizations, and it captures a conception of mathematics that is distinctly ``structural''.  These two seemingly unrelated aspects, one practical, the other philosophical, are in fact connected in a rather subtle way.  The structural character of the system, which the univalence axiom requires and indeed strengthens, permits the use of a new ``synthetic'' style of foundational axiomatics which is quite different from conventional axiomatic foundations.  One might call the conventional, set-theoretic style of foundations an ``analytic'' (or perhaps ``bottom-up'') approach, which ``analyses'' mathematical objects into constituent \emph{material} (e.g.\ sets or numbers), or at least constructs appropriate ``surrogate objects'' from such material --- think of real numbers as Dedekind cuts of rationals.  By contrast, the ``synthetic'' (or ``top-down'') approach permitted by univalent foundations is based on describing the fundamental \emph{structure} of mathematical objects in terms of their universal properties, which in type theory are given by rules of inference determining directly how the new objects map to and from all other ones.  This fundamental shift in foundational methodology has the practical effect of simplifying and shortening many proofs by taking advantage of a more axiomatic approach, as opposed to the more laborious analytic constructions.\footnote{
In a related context, it has been said that such an approach has ``all the advantages of theft over honest toil'', but the issue of how to \emph{justify} the rules for new constructions can be separated from that of their expedience.}  Indeed, in a relatively short time, a large amount of classical mathematics has already been developed in this new system:  basic homotopy theory, category theory, real analysis, the cumulative hierarchy of set theory, and many other topics.  The proofs of some very sophisticated, high-level theorems have now been fully formalized and verified by computer proof assistants --- a foundational achievement that would be very difficult to match using conventional, ``analytic'' style foundational methods.  

Indeed, this combination of a synthetic foundational methodology and a powerful computational implementation has the potential to give new life, and a new twist, to the old idea of reducing mathematics to a purely formal calculus.  Explicit formalizations that were once too tedious or complicated to be done by hand can now be accomplished in practice with a combination of synthetic methods and computer assistance.  This new formal reduction of mathematics raises again the epistemological question of whether, and in what sense, the type-theoretic basis of the formal system is purely ``logical'', and what this means about mathematics and the nature of a priori knowledge.  That is a question of significant philosophical interest, but it is perhaps better pursued independently, once the mathematical issues related to the formalization itself are more settled.

%%%%%%%%%%%%%%%%%%%%%%%%%%%%%%%%%%%%%%%%%%%%%%%%%%%%%%%%%
\section{Type theory}
%%%%%%%%%%%%%%%%%%%%%%%%%%%%%%%%%%%%%%%%%%%%%%%%%%%%%%%%%

In its current form, constructive type theory is the result of contributions made by several different people, working both independently and in collaboration.  Without wanting to give an exhaustive history (for one such, see \cite{historyoftt}), it may be said that essential early contributions were made by H.~Curry, W.~Howard, F.W.~Lawvere, P.~Martin-L\"of, D.S.~Scott, and W.W.~Tait.  

Informally, the basic system consists of the following ingredients:

\begin{itemize}
\item \myemph{Types}: $X, Y, \ldots, A\times B,\ A\rightarrow B, \ldots$, including both primitive types and type-forming operations, which construct new types from given ones, such as the product type $A\times B$ and the function type $A\rightarrow B$.
\item \myemph{Terms}: $a: A,\ b: B,\ \ldots$, including variables $x: A$ for all types, primitive terms $b:B$, and term-forming operations like $\langle a,b\rangle : A\times B$ and $\lambda x. b(x): A\rightarrow B$ associated to the type-forming operations.
\end{itemize}
One essential novelty is the use of so-called \emph{dependent types}, which are regarded as ``parametrized'' types or \emph{type families indexed over a type}.  
\begin{itemize}
\item \myemph{Dependent Types}: $x:A\vdash B(x)$  means that $B(x)$ is a type for each $x:A$, and thus it can be thought of as a function from $A$ to types.  Moreover, one can have iterated dependencies, such as:
\begin{itemize}
	\item[] $x:A\vdash B(x)$
	\item[] $x:A,\, y:B(x)\vdash C(x,y)$
	\item[] $x:A,\, y:B(x),\, z:C(x,y)\vdash D(x,y,z)$
	\item[] etc.
	\end{itemize}	
	
\item \myemph{Dependent Type Constructors}: There are special type constructors for dependent types, such as the sum $\sum_{x:A} B(x)$  and product $\prod_{x:A} B(x)$ operations.
	Associated to these are term constructors that act on dependent terms $x:A\vdash b(x):B(x)$, such as $\lambda x. b(x) : \prod_{x:A} B(x)$.
		
\item \myemph{Equations:}  As in an algebraic theory, there are then equations $s = t : A$ between terms of the same type, such as $\big(\lambda x. b(x)\big)(a) = b(a) : B(a)$.
\end{itemize}
The entire system of constructive type theory is a formal calculus of such typed terms and equations, usually presented as a deductive system by formal rules of inference.  For one modern presentation, see the appendix to \cite{hottbook}.  This style of type theory is somewhat different from the Frege-Russell style systems of which it is a descendant.  It was originally intended as a foundation for \emph{constructive} mathematics, and it has a distinctly ``predicative'' character --- for instance, it is usually regarded as open-ended with respect to the addition of new type- and term-forming operations, such as universes, so that one does not make use of the notion of ``all types'' in the way that set-theory admits statements about ``all sets'' via its first-order logical formulation.  Type theory is now used widely in the theory of programming languages and as the basis of computerized proof systems, in virtue of its good computational properties.

%%%%%%%%%%%%%%%%%%%%%%%%%%%%%%%%%%%%%%%%%%%%%%%%%%%%%%%%%
\subsection*{Propositions as types}

The system of type theory has a curious dual interpretation:
\begin{itemize}
\item On the one hand, there is the interpretation as \myemph{mathematical} objects: the types are some sort of constructive ``sets'', and the terms are the ``elements'' of these sets, which are being built up according to the stated rules of construction.
%\pause
\item But there is also a second, \myemph{logical} interpretation: the types are ``propositions'' about mathematical objects, and their terms are ``proofs'' of the corresponding propositions, which are being derived in a deductive system.
\end{itemize}
This is known as the \emph{Curry-Howard correspondence}, and it can be displayed as follows:
\medskip

\begin{center}
\begin{tabular}{|c|c|c|c|c|c|c|}
\hline 
$0$ & $1$ & $A + B$ & $A\times B^{\strut}$ & $A\rightarrow B$ & $\sum_{x:A} B(x)$ & $\prod_{x:A} B(x)$\\[1ex]
\hline 
$\bot$ & $\mathsf{T}$ & $A \vee B^{\strut}$ & $A\wedge B$ & $A\Rightarrow B$ & $\exists_{x:A} B(x)$ & $\forall_{x:A} B(x)$\\[1ex]
\hline
\end{tabular}
\end{center}

For instance, regarded as propositions, $A$ and $B$ have a conjunction $A\wedge B$, a proof of which corresponds to a pair of proofs $a$ of $A$ and $b$ of $B$ (via the $\wedge$-introduction and elimination rules), and so the terms of $A\wedge B$, regarded as a type, are just pairs $\langle a, b\rangle : A\times B$ where $a : A$ and $b : B$. Similarly, a proof of the implication $A\Rightarrow B$ is a function $f$ that, when applied to a proof $a :A$ returns a proof $f(a): B$ (\emph{modus ponens}), and so $f : A\rightarrow B$.  The interpretation of the existential quantifer $\exists_{x:A} B(x)$ mixes the two points of view:  a proof of $\exists_{x:A} B(x)$ consists of a term $a:A$ and a proof $b : B(a)$; so in particular, when it can be proved, one always has an instance $a$ of an existential statement. In classical logic, by contrast, one can use ``proof by contradiction'' to establish an existential statement without knowing an instance of it, but this is not possible here.  This gives the system has a distinctly constructive character (which can be specified in terms of certain good proof-theoretic properties).  This is one reason it is useful for computational applications.

%%%%%%%%%%%%%%%%%%%%%%%%%%%%%%%%%%%%%%%%%%%%%%%%%%%%%%%%%
\subsection*{Identity types}

Under the logical interpretation above we now have:
\begin{itemize}
\item \myemph{propositional logic}: $0$,\ $1$,\ $A + B,\ A\times B,\ A\rightarrow B$,
\item \myemph{predicate logic}: $B(x), C(x, y)$, with the \myemph{quantifiers} $\prod$ and $\sum$.
\end{itemize}
It would therefore be natural to add a primitive relation representing equality of terms $x = y$ as a \emph{type}. On the logical side, this would  represent the proposition ``x is identical to y''.  But what would it be \emph{mathematically}? How are we to continue the above table:

\begin{center}
\begin{tabular}{|c|c|c|c|c|c|c|c|}
\hline 
$0$ & $1$ & $A + B$ & $A\times B^{\strut}$ & $A\rightarrow B$ & $\sum_{x:A} B(x)$ & $\prod_{x:A} B(x)$ & $?$\\[1ex]
\hline 
$\bot$ & $\mathsf{T}$ & $A \vee B^{\strut}$ & $A\wedge B$ & $A\Rightarrow B$ & $\exists_{x:A} B(x)$ & $\forall_{x:A} B(x)$ & x = y\\[1ex]
\hline
\end{tabular}
\end{center}

We shall add to the system a new, primitive type of \emph{identity} between any terms $a,b : A$ of the same type $A$:
\[
\id{A}(a,b)\,.
\]
The \emph{mathematical} interpretation of this identity type is what leads to the homotopical interpretation of type theory.
Before we can explain that, however, we must first consider the rules for the identity types. 
 
The \myemph{introduction} rule says that $a:A$ is always identical to itself:

\begin{prooftree}
  \AxiomC{$\mathsf{r}(a):\id{A}(a,a)$}
\end{prooftree}
\medskip

%\pause

The \myemph{elimination} rule is a form of what may be called ``Lawvere's Law'':\footnote{
See \cite{Lawvere:EiHD} for a closely related principle.
}

\begin{prooftree}
 \AxiomC{$c:\id{A}(a,b)$}
  \AxiomC{$\judge{x:A}{d(x):R\bigl(x,x,\mathsf{r}(x)\bigr)}$}
  \BinaryInfC{$\mathsf{J}(a,b,c,d):R(a,b,c)$}
\end{prooftree}
\medskip
%\pause
That may look a bit forbidding when seen for the first time.  Schematically, it is saying something like:
\[
a=b\  \&\ R(x,x)\  \Rightarrow\ R(a,b)\,.
\]
Omitting the proof terms, this characterizes identity by saying that it is the \emph{least} (or better: \emph{initial}) reflexive relation.

The rules for identity types are such that if $a$\/ and $b$\/ are \emph{syntactically equal} as terms, $a=b : A$, 
then they are also \emph{identical} in the sense that there is a term $p : \id{A}(a,b)$.
But the converse is not true:  distinct terms $a\neq b$ may still be propositionally identical $p : \id{A}(a,b)$. This is a kind of \emph{intensionality} in the system, in that terms that are identified by the propositions of the system may nonetheless remain distinct syntactically, e.g.\  different polynomial expressions may determine the same function.
Allowing such syntactic distinctions to remain (rather than including a ``reflection rule'' of the form $p : \id{A}(a,b) \Rightarrow a=b$, as is done in ``extensional type theory''), gives the system its good computational and proof-theoretic properties. It also gives rise to a structure of great combinatorial complexity.

Although only the syntactically equal terms $a=b : A$ are fully interchangeable everywhere, propositionally identical ones  $p : \id{A}(a,b)$ are still interchangeable \emph{salva veritate} in the following sense: assume we are given a type family $x : A \vdash B(x)$ (regarded, if you like, as a ``predicate'' on $A$), an identity $p : \id{A}(a,b)$ in $A$, and a term  $u : B(a)$ (a ``proof of $B(a)$'').  Then consider the following derivation, using the identity rules.  %
\begin{prooftree}
     \AxiomC{$u : B(a)$}
     \AxiomC{$p : \id{A}(a,b)$}
 \AxiomC{$x : A \vdash B(x)$}
   \UnaryInfC{$x : A, y : B(x) \vdash y : B(x)$}
   \UnaryInfC{$x : A \vdash \lambda y.y : B(x) \rightarrow B(x)$}
  \BinaryInfC{$p_*: B(a) \rightarrow B(b)$}
  \BinaryInfC{$p_*u: B(b)$}
\end{prooftree}
Here $p_* = \mathsf{J}(a,b,p, \lambda y.y)$.  The resulting term $p_*u: B(b)$ (which is a derived ``proof of $B(b)$'') is called the \emph{transport} of $u$ along $p$.  Logically, this just says $$a = b\  \&\ B(a) \Rightarrow B(b)\,,$$ i.e.\ that a type family over $A$ must respect the identity relation on $A$.  As we shall see below, the homotopy interpretation provides a different view of transport; namely, it corresponds to the familiar \emph{lifting property} used in the definition of a ``fibration of spaces'':
\begin{equation}\label{diagram:fibration}
\xymatrix{
B  \ar[d] & u  \ar@{..>}[r] &  p_*u \\
A & a  \ar[r]_{p} & b
}
\end{equation}

%%%%%%%%%%%%%%%%%%%%%%%%%%%%%%%%%%%%%%%%%%%%%%%%%%%%%%%%%
\section{The homotopy interpretation}
%%%%%%%%%%%%%%%%%%%%%%%%%%%%%%%%%%%%%%%%%%%%%%%%%%%%%%%%%

Given any terms $a,b : A$, we can form the identity type  $\id{A}(a, b)$ and then consider its terms, if any, say $p, q : \id{A}(a, b)$.  Logically, $p$ and $q$ are ``proofs'' that $a$ and $b$ are identical, or more abstractly, ``reasons'' or ``evidence'' that this is so.  Can $p$ and $q$  be different?  It was once thought that such identity proofs might themselves always  be identical, in the sense that there should always be some $\alpha :  \id{\id{A}(a,b)}{(p, q)}$; however, as it turns out, this need not be so.  Indeed, there may be many distinct (i.e.\ non-identical) terms of an identity type, or none at all.  Understanding the structure of such iterated identity types is one result of the homotopical interpretation.

Suppose we have terms of ascending identity types:
\begin{align*}
  a,\ b :&\ A     \\ 
  p,\ q :&\ \id{A}(a, b)    \\ 
   \alpha,\ \beta :&\ \id{\id{A}(a,b)}{(p, q)}   \\ 
 \ldots :&\ \id{\id{\id{\ldots}}}{(\ldots)}   
\end{align*}
Then we can consider the following informal interpretation:
\begin{align*}
\text{Types} &\quad\leadsto\quad  \text{Topological spaces} \\
\text{Terms} &\quad\leadsto\quad  \text{Continuous maps} \\
a : A &\quad\leadsto\quad  \text{Points $a \in A$} \\
p : \id{A}(a, b) &\quad\leadsto\quad  \text{Paths $p$ from $a$ to $b$} \\
\alpha : \id{\id{A}(a, b)}(p, q) &\quad\leadsto\quad  \text{Homotopies $\alpha$ from $p$ to $q$} \\
& \quad\ \vdots
\end{align*}
So for instance $A$ may be a space with points $a$ and $b$, and then an identity term $p : \id{A}(a, b)$ is interpreted as a path in $A$ from $a$ to $b$, i.e.\ a continuous function $p : [0,1] \rightarrow A$ with $p0 = a$ and $p1 = b$.  If $q : \id{A}(a, b)$ is another such path from $a$ to $b$, a higher identity term $\alpha : \id{\id{A}(a, b)}(p, q)$ is then interpreted as a homotopy from $p$ to $q$, i.e.\ a ``continuous deformation'' of $p$ into $q$, described formally as a continuous function $\alpha : [0,1]\times [0,1] \rightarrow A$ with the expected behavior on the boundary of the square $[0,1]\times [0,1]$.  Higher identity terms are likewise interpreted as higher homotopies.  

Note that, depending on the choice of space $A$ and points $a, b\in A$ and paths $p, q$, it may be that there are no homotopies from  $p$ to  $q$ because, for example, those paths may go around a hole in $A$ in two different ways, so that there is no continuous way to deform one into the other.  Or there may be many different homotopies between them, for instance wrapping different numbers of times around the surface of a ball.  Depending on the space, this can become quite a complicated structure of paths, deformations, higher-dimensional deformations, etc.\ --- indeed, the investigation of this structure is what homotopy theory is all about.

One could say that the basic idea of the homotopy interpretation is just to extend the well-known topological interpretation of the \emph{simply-typed} $\lambda$-calculus \cite{AwodeyTR,AwodeyTC} (which interprets types as spaces and terms as continuous functions) to the \emph{dependently typed} $\lambda$-calculus with $\id{}$-types.  The essential new idea is then simply this: 
\begin{quote}
\emph{An identity term $p : \id{A}(a,b)$ is a path in the space~$A$ from the point $a$ to the point $b$}.
\end{quote}
Everything else essentially follows from this one idea: the dependent types $x:A \vdash B(x)$ are then forced by the rules of the type theory to be interpreted as \emph{fibrations}, in the topological sense, since one can show from the rules for identity types that the associated map $B\rightarrow A$ of spaces must have the \emph{lifting property} indicated in diagram \eqref{diagram:fibration} above (a slightly more intricate example shows that one can ``lift'' not only the endpoint, but also the entire path, and even a homotopy).  The total $\id{}$-types $\sum_{x,y:A}\id{A}(x,y)$ are naturally interpreted as \emph{path spaces} $A^I$, and the maps $f, g : A\rightarrow B$ that are identical as terms of function type $A\rightarrow B$ are just those that are \emph{homotopic} $f\sim g$. 

The homotopy interpretation was first proposed by the present author and worked out formally (with a student) in terms of Quillen model categories---a modern, axiomatic setting for abstract homotopy theory that encompasses not only the classical homotopy theory of spaces and their combinatorial models like simplicial sets, but also other, more exotic notions of homotopy (see \cite{AwodeyS:homtmi}).  The interpretation was shown to be \emph{complete} in the logical sense by Gambino and Garner \cite{gambino_garner}.\footnote{There is a technical question related to the selection of path objects  and diagonal fillers as interpretations of $\id{A}$-types and elimination $\mathsf{J}$-terms in a ``coherent'' way, i.e.\ respecting substitution of terms for variables; various solutions have been given, including \cite{warren_phd,vandenBergB:topsmi,VoevodskyV:notts,LW,Anatural}.}   These results show that intensional type theory can in a certain sense be regarded as a ``logic of homotopy'', in that the system can be faithfully represented homotopically, and then used to reason formally about spaces, continuous maps, homotopies, and so on. The next thing one might ask is, how much general homotopy theory can be expressed in this way? It turns out that a surprising amount can be captured under this interpretation, as we shall now proceed to indicate.

%%%%%%%%%%%%%%%%%%%%%%%%%%%%%%%%%%%%%%%%%%%%%%%%%%%%%%%%%
\subsection*{The fundamental groupoid of a type}

Like path spaces in topology, identity types endow each type with the structure of a \emph{groupoid}: a category in which every arrow has an inverse.  
$$
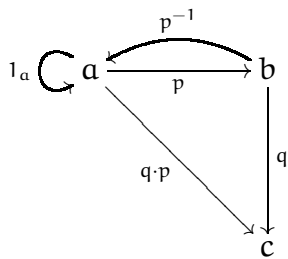
\xymatrix{ 
a \ar@(ul,dl)_{1_a} \ar[rr]_{p}  \ar[ddrr]_{q\cdot p} && b  \ar@/{}_{1pc}/[ll]_{p^{-1}}  \ar[dd]^{q} \\
\\
&& c} 
$$

The familiar \emph{laws of identity}, namely reflexivity, symmetry, and transitivity are provable in type theory, and their proof terms therefore act on identity terms, providing the \emph{groupoid operations} of unit, inverse, and composition:
\begin{align*}
r : \id{}(a,a) &\qquad \text{reflexivity} && a \to a \\
s: \id{}(a,b)\rightarrow\id{}(b,a) &\qquad \text{symmetry} && \xymatrix{ a \ar[r] & b \ar@/{}_{.5pc}/[l]} \\
t :  \id{}(a,b)\times\id{}(b,c)\rightarrow\id{}(a,c) &\qquad \text{transitivity} && \xymatrix{ 
	a \ar[r]  \ar[dr] & b \ar[d] \\
	& c} 
\end{align*}
The \emph{groupoid laws} of units, inverses, and associativity also hold ``up to homotopy'', i.e.\ up to the existence of a higher identity term.
This means that instead of e.g.\ $p^{-1}\cdot p = 1_a$, we have a higher identity term:
\[
\alpha :  \id{\Id}\left(p^{-1}\cdot p, 1_a\, \right)
\]
as indicated in:
$$
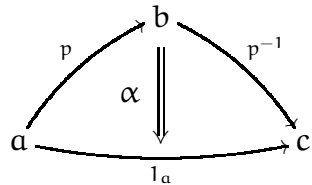
\begin{array}{c}
\ \xy
(0,0)*+{a}="a";
(450,400)*+{b}="b";
(900,0)*+{c}="c";
{\ar@/^.5pc/^{p} "a";"b"};
{\ar@/^.5pc/^{p^{-1}} "b";"c"};
{\ar@/_.5pc/_{1_a} "a";"c"};
{\ar@{=>} (450,300)*{};(450,0)*{}};(350,150)*{\alpha};
\endxy \ 
\end{array}
$$
Indeed, this is just the same situation that one encounters in defining the fundamental group of a space in classical homotopy theory, where one shows e.g.\ that composition of paths is associative up to homotopy by reparametrization of the composites.  In fact, in virtue of the homotopy interpretation, the classical case is really just an instance of the more general, type theoretic one.

Inspired by this occurence of type theoretic groupoids, Hofmann and Streicher \cite{HofmannM:gromtt} discovered an interpretation of the entire system of type theory into the category of all groupoids, which was a precursor of the homotopy interpretation.  It was used, for instance, to estable the above mentioned fact that identity types may have elements that are not themselves identical.

The identity structure of a general type may actually be much richer than that of just a groupoid;  as in homotopy theory, there may be non-trivial higher identities, representing higher homotopies between homotopies, and this structure may go on to higher and higher identities without ever becoming degenerate.
 
 $$
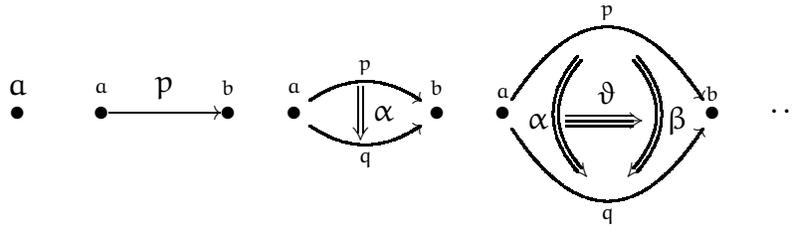
\begin{array}{c}
\begin{array}{ccccc}
\ \xy
(0,0)*{\bullet};
(0,80)*{a};
\endxy \quad
&
\ \xy
(0,0)*{\bullet}="a";
(0,80)*{\scriptstyle a};
(400,0)*{\bullet}="b";
(400,80)*{\scriptstyle b};
{\ar "a";"b"};
(200,80)*{p};
\endxy \ 
&
\ \xy
(0,0)*+{\bullet}="a";
(0,80)*{\scriptstyle a};
(450,0)*+{\bullet}="b";
(450,80)*{\scriptstyle b};
{\ar@/^1pc/^{p} "a";"b"};
{\ar@/_1pc/_{q} "a";"b"};
{\ar@{=>} (210,85)*{};(210,-85)*{}};
(280,0)*{\alpha};
\endxy \ 
&
\xy 0;/r.22pc/: 
(0,15)*{}; 
(0,-15)*{}; 
(0,8)*{}="A"; 
(0,-8)*{}="B"; 
%{\ar@{=>} "A" ; "B"}; 
{\ar@{=>}@/_.75pc/ "A"+(-4,1) ; "B"+(-3,0)}; 
(-10,0)*{\alpha};
{\ar@{=}@/_.75pc/ "A"+(-4,1) ; "B"+(-4,1)}; 
{\ar@{=>}@/^.75pc/ "A"+(4,1) ; "B"+(3,0)}; 
(10,0)*{\beta};
{\ar@{=}@/^.75pc/ "A"+(4,1) ; "B"+(4,1)}; 
{\ar@3{->} (-6,0)*{} ; (6,0)*+{}}; 
(0,4)*{\vartheta};
%{\ar@3{->} (2,0)*{} ; (6,0)*+{}}; 
(-15,1)*+{\bullet}="1"; 
(-15,4)*{\scriptstyle a};
(15,1)*+{\bullet}="2"; 
(15,4)*{\scriptstyle b};
{\ar@/^2.75pc/^{p} "1";"2"}; 
{\ar@/_2.75pc/_{q} "1";"2"}; 
\endxy \
&
\ \xy
\dots
\endxy \ 
\end{array} 
\end{array}
$$
The resulting structure is that of an \emph{$\omega$-groupoid}, which is something that  has also appeared  elsewhere in mathematics --- twice!  As already mentioned, such ``infinite-dimensional groupoids'' also occur in homotopy theory, where the fundamental $\omega$-groupoid of a space is an algebraic  invariant that respects the homotopy type (according to Grothendieck's famous ``homotopy hypothesis'' these groupoids contain all the essential information of the space up to homotopy); but also in category theory, one has considered the idea of an \emph{$\omega$-category}, with not only objects and arrows between them, but also 2-arrows between arrows, 3-arrows between 2-arrows, and so on.  It is indeed remarkable that the same notion has now appeared again in logic, as exactly the structure of iterated identity in type theory.\footnote{See \cite{LumsdaineP:weaci} and \cite{vandenBergB:typwg} for  details.}

%%%%%%%%%%%%%%%%%%%%%%%%%%%%%%%%%%%%%%%%%%%%%%%%%%%%%%%%%%%%%%%%
\subsection*{Homotopy levels}

One of the most useful new discoveries is that the system of all types is naturally stratified into ``homotopy levels'' by a hierarchy of definable conditions.\footnote{This concept is due to Voevodsky, cf.~\cite{VoevodskyV:notts}.  Also see \cite{hottbook}, ch.~7.}
 
 At the lowest level are those types that are \emph{contractible} in the following sense.
  \[
  \text{$X$ is contractible $=_\text{def}$} \quad \sum_{x:X}\prod_{y:X}\id{X}(x,y)\,.
  \]
 Under the logical reading, this condition says that $X$ is a ``singleton'', in that there is an element $x:X$ such that everything $y:X$ is identical to it.  So roughly, these are the types that have just one element, up to homotopy.
 
 The next level consists of the \emph{propositions}, defined as those types whose identity types are always contractible,
 \[
 \text{$X$ is a proposition $=_\text{def}$} \quad \prod_{x,y:X}\ \mathsf{Contr}(\id{X}(x,y))\,.
 \]
It is not hard to see that such types are contractible \emph{if} they are inhabited---thus they are like ``truth values'', either false (i.e.\ empty) or true (i.e.\ contractible), and then essentially uniquely so.  In other words, the elements of a proposition contain no further information, other than the mere inhabitation of the proposition, which we interpret to mean that it holds.

At the next level are the \emph{sets}, which are types whose identity relation is always a proposition:
\[
\text{$X$ is a set $=_\text{def}$} \quad \prod_{x,y:X}\ \mathsf{Prop}(\id{X}(x,y))
\]
These types have the familiar, set-like behavior that the identity proofs, when they exist, are unique (again ``up to homotopy''). 

Next come the types whose identity types are sets, which may be called \emph{groupoids}, because they are like the algebraic groupoids just discussed:
\[
\text{$X$ is a groupoid $=_\text{def}$} \quad \prod_{x,y:X}\ \mathsf{Set}(\id{X}(x,y))
\]
These types may have distinct identity proofs between elements, but all higher identity proofs are degenerate.  

The general  pattern is now clear:
\[
\text{$X$ has homotopy level $n+1$ $=_\text{def}$} \quad \prod_{x,y:X}\ \mathsf{H}_{n}(\id{X}(x,y))
\]
Thus the  types $X$ of homotopy level $n+1$ (for which we write $\mathsf{H}_{n+1}(X)$) are the types whose identity relation is of homotopy level $n$; these types correspond to the higher-dimensional groupoids of category theory, when we think of identity terms as higher-dimensional arrows.  To start the numbering, we may set the contractible types to be level 0.

The homotopy level of a type is the height at which the tower of iterated identity types becomes degenerate; under the homotopy interpretation this corresponds (up to a shift in numbering) to the notion of a space being a \emph{homotopy $n$-type}, which is usually defined as the greatest $n$ such that the $n$-th homotopy group is non-trivial.  In each case, it is a measure of the complexity of the type/space --- in the former case in terms of higher identities, and in the latter in terms of higher homotopies.

The recognition that types have these different degrees of complexity allows for a more refined version of the propositions-as-types idea, according to which only those types that are ``propositions'' in the sense of the homotopy levels are read as bare \emph{assertions}, while others are regarded more discriminately as \emph{structured objects} of various kinds.  Accordingly, a type family $x:A \vdash B(x)$ such that all values $B(x)$ are propositions can be regarded as a simple ``predicate'' (or a ``relation'' depending on the arity), while a family of sets, groupoids, etc.\ is viewed more accurately as a structure on $A$. 

The stratification of types by homotopy levels gives us a new view of the mathematical universe, which is now seen to be  arranged not only into the familiar, one-dimensional hierarchy of \emph{size}, determined by a system of universes $\mathcal{U}_0, \mathcal{U}_1, \mathcal{U}_2, \dots$, but also into a hierarchy of \emph{homotopy levels}, which form a second dimension independent of the first (see Fig.~1).

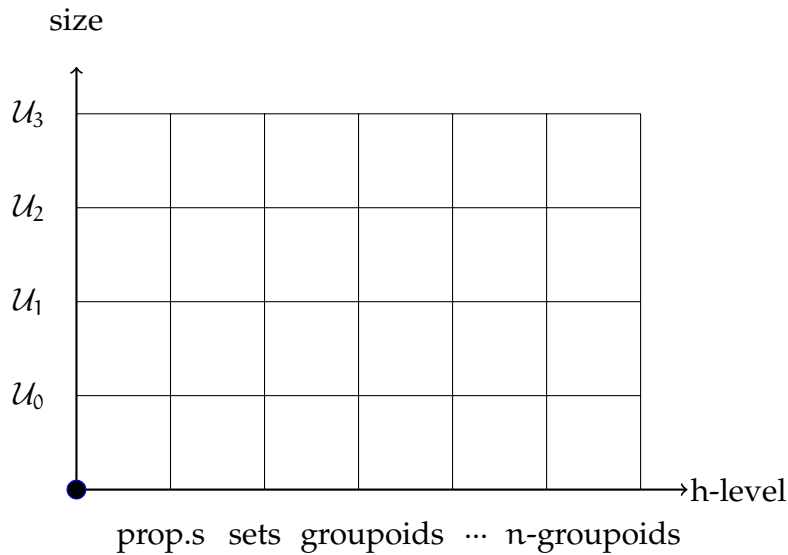
\begin{figure}
\begin{center}
\begin{tikzpicture}[yscale=1.25,xscale=1.25] 
%\newcounter{sub}
%\setcounter{sub}{1}
\draw (0, 0) grid (6, 4);  % 5x5 grid
%%  for each loop
%\foreach \y in {3,2,1}{
%\foreach \x in {1,2,3}{
%\draw [color=blue, fill=blue] (\x, \y) circle (0.1);  % for the blue cicles
%\node at (\x+0.25, \y+0.25) {$P_{\thesub}$};          % for the P_sub labels
%\stepcounter{sub}
%}}                                                    % what follows are OP's code mainly
\node at (.9,-.55) {prop.s};
\node at (1.9,-.48) {sets};
\node at (3.15,-.5) {groupoids};
\node at (4.25,-.5) {...};
\node at (5.5,-.5) {$n$-groupoids};
\node at (-.5,1) {$\mathcal{U}_0$};
\node at (-.5,2) {$\mathcal{U}_1$};
\node at (-.5,3) {$\mathcal{U}_2$};
\node at (-.5,4) {$\mathcal{U}_3$};
%
%\node at (-0.1, -0.5) {};% origin label
%\draw [green, thick, domain=0:4] plot ({0},\x);
%\draw [green, thick, domain=0:4] plot ({4},\x);

\node at (7, 0) {\ h-level};% x-axis label
\draw [thick,->] (0, 0) -- (6.5, 0);
\node at (0, 5) {size};% y-axis label
\draw [thick,->] (0, 0) -- (0, 4.5);
\draw [color=blue, fill=black] (0, 0) circle (0.1);
\end{tikzpicture}
\caption{The 2D hierarchy of types}
\end{center}
\end{figure}

%%%%%%%%%%%%%%%%%%%%%%%%%%%%%%%%%%%%%%%%%%%%%%%%%%%%%%%%%%%%%%%%
\section{Higher inductive types}
%%%%%%%%%%%%%%%%%%%%%%%%%%%%%%%%%%%%%%%%%%%%%%%%%%%%%%%%%

The recognition and use of the notion of homotopy level of a type has made the entire system of type theory more expressive and powerful, for example by allowing greater control over the introduction of new type constructions.   One such construction that was formerly problematic but is now better behaved is the construction of the \emph{quotient type} $A/\!\sim$ of a type $A$ by an equivalence relation $x, y:A\vdash x\sim y$.  When $x\sim y$ is known to be a proposition for all $x, y:A$, then the quotient $A/\!\sim$ will be a set, and the introduction and elimination rules can be determined without difficulty.  Such ``set quotients'' can be constructed, roughly speaking, as equivalence classes \cite{VoevodskyV:notts}; or they can be introduced axiomatically \cite{hottbook}, essentially by stating rules that say that the identity type of $A/\!\sim$ is a \emph{relation} (i.e.\ a family of propositions) that is freely generated by the equivalence relation $x\sim y$.

The latter, axiomatic approach is a special case of the very powerful construction method of \emph{higher inductive types}, which are a systematic way of introducing new types with stipulated points, paths, higher paths, etc..  In order to explain this further, let us first recall how type theory deals with \emph{ordinary} inductive types, like the natural numbers.   The natural numbers $\mathbb{N}$ can be implemented as an inductive type via rules that may be represented schematically as:
\[
			\mathbb{N} := \begin{cases} &0 : \mathbb{N}\\
		 				&s : \mathbb{N} \rightarrow \mathbb{N}
						\end{cases}
\]
The terms $0$ and $s$ are the \emph{introduction rules}  for this type.
The recursion property of $\mathbb{N}$ is captured by an \emph{elimination rule}:
\begin{prooftree}
% \AxiomC{$X$ type}
 \AxiomC{$a : X$}
  \AxiomC{$f : X \rightarrow X$}
%   \UnaryInfC{$x : A, y : B(x) \vdash y : B(x)$}
%   \UnaryInfC{$x : A \vdash \lambda y.y : B(x) \rightarrow B(x)$}
%  \BinaryInfC{$p_*: B(a) \rightarrow B(b)$}
  \BinaryInfC{$\mathsf{rec}(a,f) :  \mathbb{N} \rightarrow X$}
\end{prooftree}
which says that given any structure of the same kind as $\mathbb{N}$, there is a map $\mathsf{rec}(a,f)$ to it from $\mathbb{N}$, which furthermore preserves the structure, as stated by the following \emph{computation rules}:
\begin{align*}
\mathsf{rec}(a,f)(0) &= a\,,\\
\mathsf{rec}(a,f)(sn) &= f(\mathsf{rec}(a,f)(n))\,.
\end{align*}

The map $\mathsf{rec}(a,f) :  \mathbb{N} \rightarrow X$ is actually required to be the \emph{unique} one satisfying the computation rules, a condition that can be ensured either with a further computation rule or by reformulating the elimination rule as a more general \emph{induction principle} rather than a recursion principle (cf.~\cite{wtypes}).

In more algebraic terms, one would say that $(\mathbb{N}, 0, s)$ is the \emph{free structure} of this kind.
We remark that it can be shown on the basis of these rules, and without further assumptions, that \emph{$\mathbb{N}$ is a set} in the sense of the hierarchy of homotopy levels.

%%%%%%%%%%%%%%%%%%%%%%%%%%%%%%%%%%%%%%%%%%%%%%%%%%%%%%%%%
\subsection*{The circle}

We now want to use the same method of specifying a new type by introduction and elimination rules (which amount to specifying the mappings to and from other types), but now with generating data that may include also elements of identity types, in addition to elements of the type itself and operations on it.  A simple example is the following.

The homotopical circle $\Sn$ can be given as an inductive type involving one ``base point'' and one ``higher-dimensional'' generator:
\[
			\Sn := \begin{cases} &\base : \Sn\\
					&\lloop : \id{\Sn}(\base,\base)
					\end{cases}
\]
The element $\lloop :\id{\mathbb{S}}(\base,\base)$ can therefore be regarded as a ``loop'' at the basepoint $\base: \Sn$, i.e.\ a path that starts and ends at $\base$.  
The corresponding recursion property of $\Sn$ is then given by the following elimination rule,
\begin{prooftree}
% \AxiomC{$X$ type}
 \AxiomC{$a : X$}
  \AxiomC{$p : \id{X}(a,a)$}
%   \UnaryInfC{$x : A, y : B(x) \vdash y : B(x)$}
%   \UnaryInfC{$x : A \vdash \lambda y.y : B(x) \rightarrow B(x)$}
%  \BinaryInfC{$p_*: B(a) \rightarrow B(b)$}
  \BinaryInfC{$\mathsf{rec}(a,p) :  \mathbb{S} \rightarrow X$}
\end{prooftree}
with computation rules,
\begin{align*}
\mathsf{rec}(a,p)(\base) &= a\,,\\
\mathsf{rec}(a,p)_!(\lloop) &= p\,.
\end{align*}
There is an obvious analogy to the rules for $\mathbb{N}$.\footnote{
A map $f : A\rightarrow B$ induces a map  on identities, taking each $p:\id{A}(a,b)$ to a term in $\id{A}(fa,fb)$ which we here write $f_!p$ (see \cite{hottbook}, ch.~2).
}  
The map $\mathsf{rec}(a,p) :  \mathbb{S} \rightarrow X$ is then moreover required to be unique up to homotopy, which again is achieved either with additional computation rules or a generalized elimination rule in the form of ``circle induction'' rather than ``circle recursion'' (see \cite{hottbook},\cite{Sojakovatech_report}).

Conceptually, these rules suffice to make the structure $(\mathbb{S}, \base, \lloop)$ into the ``free type with a point and a loop''. 
To see that it actually behaves as it should to be the homotopical circle, one can verify that it has the correct homotopy groups (cf.~\cite{licata_shulman}):

\begin{theorem}[Shulman 2011] The type-theoretic circle $\Sn$ has the following homotopy groups:
\[
\pi_n(\mathbb{S}) = \begin{cases}
	\mathbb{Z}, &\text{if $n = 1$,}\\
	0, &\text{if $n\neq 1$.}
	\end{cases}
\]
\end{theorem}
The homotopy groups $\pi_n(X,x)$ for any type $X$ and basepoint $x:X$ can be defined as usual in terms of loops at $x$ in $X$, i.e.\  identity elements $\id{X}(x,x)$, ``modulo homotopy'', i.e.\ modulo higher identities.  The proof of the above theorem can be given entirely within the system of type theory, and it combines methods from classical homotopy theory with ones from constructive type theory in a novel way, using Voevodsky's univalence axiom (a sketch is given in Section \ref{synthetic} below).  The entire development has been fully formalized \cite{licata_shulman}.

%%%%%%%%%%%%%%%%%%%%%%%%%%%%%%%%%%%%%%%%%%%%%%%%%%%%%%%%%
\subsection*{The interval}

The homotopical interval $\mathbb{I}$ is also a higher inductive type, this time generated by the basic data:
\[
			\mathbb{I} :=  \begin{cases} 
						&0,1 : \mathbb{I}\\
						&\mathsf{path}: \id{\mathbb{I} }(0,1)
						\end{cases}
\]
Thus $\mathsf{path}:\id{\mathbb{I}}(0,1)$ represents a path from $0$ to $1$ in $\mathbb{I}$. The elimination and computation rules are analogous to those for the circle, but now with separate endpoints $0$ and $1$.  So given any path $p:\id{X}(a,b)$ between points $a$ and $b$ in any type $X$, there is a unique (up to homotopy) map $\mathbb{I} \rightarrow X$ taking $0$ to $a$, $1$ to $b$, and $\mathsf{path}$ to $p$.
This specification makes the structure $(\mathbb{I}, 0, 1, \mathsf{path})$ the ``free type with a path''.

In terms of this example, we can plainly compare the methodology behind the use of higher inductive types in homotopy type theory with the conventional approach of classical topology:
\begin{quote}
In classical topology, we start with the \emph{interval} and use it to define the notion of a \emph{path}.\\
In homotopy type theory, we start with the notion of a \emph{path}, and use it to define the \emph{interval}.
\end{quote}
The notion of a \emph{path}, recall, is a primitive one in our system, namely a term of identity type.  In terms of these, one can then determine the interval $\mathbb{I}$ via its mappings, rather than the other way around.

%%%%%%%%%%%%%%%%%%%%%%%%%%%%%%%%%%%%%%%%%%%%%%%%%%%%%%%%%
\subsection*{Constructing higher inductive types}

The higher inductive types mentioned so far were introduced \emph{axiomatically}, by stating their basic rules.  
Can we instead \emph{construct} them, similarly to the way that quotients can also be constructed from equivalence classes?  One possible way to do this is by what is sometimes called an ``impredicative encoding'', which is a construction that involves a quantification over ``all types''. 

Consider first some related examples that are not \emph{higher} inductive types, but which can be determined by such impredicative encodings.  First let $p$ and $q$ be propositions and consider:
\begin{align*}
p \vee q \ &=_{\mathrm{def}} \ \forall_{x}\big[ (p \imp x)\wedge(q\imp x)\imp x \big]
%A + B \ &= \ \prod_{X}\big[(A \arr X)\times(B\arr X)\arr X\big]
\end{align*}
where the quantifier $\forall_{x}$ is over \emph{all propositions} $x$ (in a given universe). Among propositions, this type has the correct behavior to be the disjunction of $p$ and $q$.

Next let $A$ and $B$ be sets, and consider: 
%\pause
\begin{align*}
%p \vee q \ &= \ \forall_{x}\big[ (p \imp x)\wedge(q\imp x)\imp x \big]\\
A + B \ &=_{\mathrm{def}}\ \prod_{X}\big[(A \arr X)\times(B\arr X)\arr X\big]
\end{align*}
where the product $\prod_{X}$ is over \emph{all sets} $X$ (again, in a given universe).  In order for this to actually be the coproduct among all sets, this specification requires an additional \emph{coherence condition} saying that the transformations 
$$\alpha_X : \big((A \arr X)\times(B\arr X)\big)\ \to\ X$$
are \emph{natural in $X$}, in a straightforward sense that we will not spell out here.

The same general idea can be used for the interval and the circle: for a type $X$, define the \emph{path space} $I(X)$ and the \emph{(free) loop space} $L(X)$ by:
\begin{align*}
I(X)\ &=_{\mathrm{def}}\ \sum_{x,y:X}\id{X}(x,y) \\
L(X)\ &=_{\mathrm{def}}\ \sum_{x:X}\id{X}(x,x)
\end{align*}
Morally, we expect that $I(X) = \mathbb{I}\arr X\ $\ and\ $\ L(X) = \Sn\arr X$, so as a first approximation we set:
\begin{align*}
\mathbb{I} \ &=_{\mathrm{def}}\ \prod_{X}\big[I(X)\arr X\big]\\[1ex]
\Sn \ &=_{\mathrm{def}}\ \prod_{X}\big[L(X)\arr X\big]
\end{align*}
%\pause
where the product $\prod_{X}$ is now over all \emph{groupoids} $X$ (types of homotopy level 3, in a given universe).  In order to get the correct elimination rules for these types, we again add a further \emph{coherence condition}, now involving higher-order naturality, which again we will not spell out here.\footnote{
See \cite{Awodey:Impredicative}.  The displayed formula for the circle was first considered by M.~Shulman.}

The possibility of a ``logical construction of the circle'' and similar constructions of some other higher inductive types are current work in progress.  At present they require either a general assumption of ``impredicativity'', or more specialized ``resizing rules'', or some other device to handle the shift in universes involved in the quantification over ``all types''.

\bigskip

Many basic spaces and constructions can be introduced directly as higher inductive types.  These include, for example:

\begin{itemize}
\item higher spheres $\mathbb{S}^n$, mapping cylinders, tori, cell complexes,
\item suspensions $\Sigma A$, homotopy pushouts,
\item truncations, such as connected components $\pi_0(A)$ and  ``bracket'' types $[A]$ (cf.~\cite{AwodeyS:prot}),
\item (higher) homotopy groups $\pi_n$, Eilenberg-MacLane spaces $K(G,n)$, Postnikov systems,
\item a Quillen model structure on the system of all types,
\item quotients by equivalence relations and more general quotients,  
\item free algebras, algebras presented by generators an relations, 
\item the real numbers, the surreal numbers,
\item the cumulative hierarchy of Zermel-Fraenkel sets.
\end{itemize}

The use of higher inductive types is a topic that is curently under very active investigation (see e.g.~\cite{LumsdaineP:higit}).

%%%%%%%%%%%%%%%%%%%%%%%%%%%%%%%%%%%%%%%%%%%%%%%%%%%%%%%%%%%%%%%
\section{Univalence}
%%%%%%%%%%%%%%%%%%%%%%%%%%%%%%%%%%%%%%%%%%%%%%%%%%%%%%%%%

Voevodsky has proposed a new foundational axiom to be added to type theory: the \emph{univalence axiom}.
It is motivated by the homotopy interpretation and makes precise the informal mathematical practice of ``identifying'' isomorphic objects.  
Especially when combined with higher inductive types, this new axiom is a powerful addition to the system. Although it is formally incompatible with the naive interpretation of type theory according to which all types are sets, it is provably consistent with the homotopical interpretation \cite{ssets}.  Its status as a constructive principle is still unsettled, however, and that question is the focus of much current research.

\subsection*{Isomorphism, equivalence, and invariance}

In type theory, the notion of a \emph{type isomorphism} $A \cong B$ is definable as usual; namely, the statement
\[
\text{there are}\ f:A\toarrow B\ \text{and}\ g: B\toarrow A\ \text{with}\ gf(x) = x\  \text{and}\ fg(y) = y
\]
is formalized by the type of isomorphisms,
\[
\mathrm{Iso}(A,B)\ =_{\mathrm{def}}\  \sum_{f:A\toarrow B}\sum_{g: B\toarrow A}\big(\prod_{x:A}\id{A}(gf(x),x)\times\prod_{y:B}\id{B}(fg(y),y)\big)\,.
\]
Under the logical reading, this type expresses exactly the preceding informal statement. The types $A$ and $B$ are isomorphic just if this type is inhabited by a term, which is then exactly an isomorphism between $A$ and $B$.  Here we see the propositions-as-types idea at work: a proof of the proposition $A\cong B$ is the same thing as a term of the type  $\mathrm{Iso}(A,B)$, namely, an isomorphism.

There is also a more refined notion of \emph{equivalence} of types
$
A \simeq B
$
which adds a further  ``coherence'' condition relating the identity terms of $\id{A}(gf(x),x)$ and $\id{B}(fg(y),y)$ via $f$ and $g$ (see \cite{hottbook}, chapter 4).  Since every isomorphism can be ``promoted'' to an equivalence, the latter condition is no ``stronger'' logically; nonetheless, it is worth the extra trouble to consider, because being an equivalence $f:A\simeq B$ is always a propositional condition, whereas being an isomorphism $f:\mathrm{Iso}(A,B)$ need not be one.  Under the homotopy interpretation, the type $A\simeq B$ consists of the \emph{homotopy equivalences} of spaces. 
The notion of type equivalence also subsumes \emph{categorical equivalence} (for groupoids), \emph{isomorphism} (for sets), and \emph{logical equivalence} (for propositions).

Now, it is an important fact about type theory that all ``definable properties'' $P(X)$ of types $X$ (formally, any type expression with a type variable $X$) can be shown to respect type equivalence, in the sense that  $A \simeq B$ and $P(A)$ imply $P(B)$; indeed, if $A \simeq B$ then $P(A) \simeq P(B)$.  Briefly, we may say that all type-theoretic properties and concepts are \emph{invariant}.\footnote{
This of course does not hold in set theory.  For example, consider the sets $\{\emptyset\}$ and $\{\{\emptyset\}\}$, which are isomorphic but are distinguished by the property $P(X) = (\exists x,y)\,x\in y\in X$.
}  
It therefore follows that equivalent types $A \simeq B$ are \emph{indiscernable} within the system.  Thus it is natural to ask how equivalence is related to the \emph{identity} of the types $A$ and $B$.

%%%%%%%%%%%%%%%%%%%%%%%%%%%%%%%%%%%%%%%%%%%%%%%%%%%%%%%%%
\subsection*{The univalence axiom}

To reason internally about identity of types $A$ and $B$, we need to add to the basic system a \emph{type universe}\ $\U$, with an identity type,
\[
\id{\U}(A,B)
\]
expressing the relation of identity of types.  The usual rules for identity then imply that identity of types implies their equivalence (because equivalence is reflexive), and so there is a comparison map,
\[
\id{\U}(A,B) \toarrow (A \simeq B).
\]

The \emph{univalence axiom} asserts that this map is itself an equivalence:
\begin{equation}\tag{UA}
\id{\U}(A,B)\ \simeq\ (A \simeq B)
\end{equation}
%\pause
%
So UA can be read: \emph{identity is equivalent to equivalence.}  It internally identifies those types that are equivalent, and therefore indiscernable.  Indeed, since UA is an equivalence, there is a map coming back:
\begin{equation*}
 \id{\U}(A,B) \longleftarrow  (A\simeq B)\,.
\end{equation*}
Regarded as a map of types, we write this as $\mathsf{ua} : (A\simeq B) \rightarrow \id{\U}(A,B)$, which maps equivalences of types to identities between them.  Read logically, $\mathsf{ua}$ is a proof of the statement that \emph{equivalent types are identical}; thus in particular, isomorphic sets, groups, etc., are also identified.
%\pause

Note that in the extended system with a universe $\U$, the univalence axiom is just what is needed to maintain the above-mentioned \emph{invariance} of all ``properties''~$P(X)$:
\[
A \simeq B\ \text{implies}\ P(A) \simeq P(B)\,,
\]
for we can take $P(X) = \id{\U}(A,X)$ to see that equivalent types must be identified.

Note also that UA implies that $\U$, in particular, is not a set: for there are two distinct isomorphisms $2\cong 2$, and these therefore correspond by UA to two distinct identity terms in $\id{\U}(2,2)$.

Finally, we mention that the computational character of UA is still an open question. The system of type theory without it has some desirable properties, like the so-called ``strong normalization'' of terms, which implies the decidability of the syntactic equality relation of terms $a = b :A$.  Adding a new axiom like univalence is likely to disrupt this property, but does it completely destroy the constructive character of the system?  This is one of the open questions currently under active investigation.

%%%%%%%%%%%%%%%%%%%%%%%%%%%%%%%%%%%%%%%%%%%%%%%%%%%%%%%%%%
\section{Synthetic reasoning}\label{synthetic}
%%%%%%%%%%%%%%%%%%%%%%%%%%%%%%%%%%%%%%%%%%%%%%%%%%%%%%%%%%

In homotopy type theory, what is called the ``synthetic'' style of reasoning involves making use of the \emph{primitive geometric element} introduced by taking the notion of a \emph{path} as basic, rather than reducing it to maps from the real interval $[0,1]$.  This method is especially powerful in combination with the univalence axiom.  By way of example, let us sketch the proof of the above-mentioned theorem that the  fundamental group of the circle $\Sn$ is the integers $\Z$.

%%%%%%%%%%%%%%%%%%%%%%%%%%%%%%%%%%%%%%%%%%%%%%%%%%%%%%%%%
\subsection*{Computing $\pi_1\Sn$}

\begin{figure}\centering
  \begin{tikzpicture}[xscale=1.4,yscale=.6]
    \node (R) at (2,1) {$\mathbb{R}$};
    \node (S1) at (2,-2) {$\Sn$};
    \draw[->] (R) -- node[auto] {$\cov$} (S1);
    \draw (0,-2) ellipse (1 and .4);
    \draw[dotted] (1,0) arc (0:-30:1 and .8);
    \draw (1,0) arc (0:90:1 and .8) arc (90:270:1 and .3) coordinate (t1);
    \draw[white,line width=4pt] (t1) arc (-90:90:1 and .8);
    \draw (t1) arc (-90:90:1 and .8) arc (90:270:1 and .3) coordinate (t2);
    \draw[white,line width=4pt] (t2) arc (-90:90:1 and .8);
    \draw (t2) arc (-90:90:1 and .8) arc (90:270:1 and .3) coordinate (t3);
    \draw[white,line width=4pt] (t3) arc (-90:90:1 and .8);
    \draw (t3) arc (-90:-30:1 and .8) coordinate (t4);
    \draw[dotted] (t4) arc (-30:0:1 and .8);
    \node[fill,circle,inner sep=1pt,label={below:\scriptsize $\base$}] at (0,-2.4) {};
    \node[fill,circle,inner sep=1pt,label={above left:\scriptsize 0}] at (0,.2) {};
    \node[fill,circle,inner sep=1pt,label={above left:\scriptsize 1}] at (0,1.2) {};
    \node[fill,circle,inner sep=1pt,label={above left:\scriptsize 2}] at (0,2.2) {};
  \end{tikzpicture}
  \caption{The winding map in classical topology}\label{fig:winding}
\end{figure}
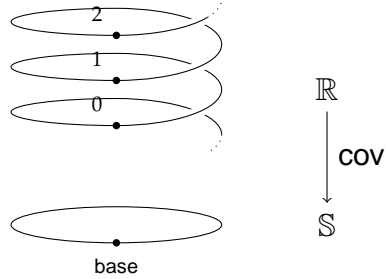

To compute the fundamental group of the circle $\Sn$, just as in the classical proof, we shall make use of the ``universal cover'' (see Fig.~\ref{fig:winding}).
%\pause
%
%\begin{definition}[Universal Cover of $\Sn^1$]
As a covering space and therefore a fibration, the universal cover will be a dependent type over $\Sn$, which, in the presence of a universe $\U$, is simply a map $$\cov : \Sn \to \U.$$  We can define such a type family using the recursion property of the circle; indeed, we just need the following data:
\begin{itemize}
\item a point $A:\U$
\item a loop $p : \id{\U}(A,A)$
\end{itemize}
For the point $A$ we shall take the integers $\Z$.
By univalence, to give a loop $p : \id{\U}(A,A)$ in $\U$, it suffices to give an equivalence $q : \Z\simeq \Z$.
But since $\Z$ is a set, equivalences are just isomorphisms; so for $q$ we can take the successor function
$\suc : \Z\cong\Z$.

\begin{definition*}[Universal cover of $\Sn$]
The  type family $\cov : \Sn \to \U$ is given by circle recursion with 
\begin{align*}
    \cov(\base) &= \Z\,,\\
    \cov(\lloop) &= \ua(\suc)\,.
\end{align*}
\end{definition*}

As in classical homotopy theory, we then use the universal cover to define the ``winding number'' of any path $p : \id{\Sn}(\base,\base)$ by $\mathsf{wind}(p)=p_*(0)$, where $p_*$ is the transport operation along $p$.
 This gives a map from the type $\Omega(\Sn)$ of (based) loops in $\Sn$ to the integers,
\[
\mathsf{wind} : \Omega(\Sn)\to\Z,
\]
which can be shown to be inverse to the map $\Z\to\Omega(\Sn)$ defined by composing $\lloop$ with itself a given number  of times $i$,
\[
i \mapsto \lloop^i.
\]

This proof can be formalized in a very efficient way, and the result is no longer than a conventional, ``unformal'' mathematical proof (see \cite{licata_shulman}).  This is a real advance over the traditional ``analytic'' style of formalization, which would require defining the circle as a subspace of the Euclidean plane $\mathbb{R}^2$, defining homotopies via continuous maps from the unit interval $[0,1]$, using reparametrizations of paths to define their composition, defining the real numbers $\mathbb{R}$ and the winding map via trigonometric functions, and so on.  To be sure, those are worthwhile mathematical objects and constructions in their own right!  But by avoiding the need for them in this case, the synthetic approach seems to be somehow closer to the real ``essence'' of the homotopical fact being proved.

%%%%%%%%%%%%%%%%%%%%%%%%%%%%%%%%%%%%%%%%%%%%%%%%%%%%%%%%%
\subsection*{The cumulative hierarchy of sets}

As a final example of synthetic reasoning in the full system of homotopy type theory with higher inductive types and univalence, we consider a somewhat ``experimental'' construction which gives the cumulative hierarchy of sets; see \cite{hottbook} for the details.   

Given a universe $\U$, we  make the \emph{cumulative hierarchy} $V$ of sets in $\U$ as a higher inductive type with the following generating constructors (we shall write $x=_V y$ for $\id{V}(x,y)$ for easier comparison with more familiar treatments):
\begin{enumerate}
\item For any  type $A:\U$ and any map $f : A\arr V$, there is a ``set'',
\[
\mathsf{set}(A,f) : V\,.
\]
We think of $\mathsf{set}(A,f)$ as the image of $A$ under $f$, i.e.\ the classical set $$\{ f(a)\ |\ a \in A \}.$$

\item \label{bisim} For all $A : \U$ and $f : A \arr V$ and $B : \U$ and $g : B \arr V$ such that
    \[
      \big(\forall{a:A}\,\exists{b:B}\ f(a) =_V g(b)\big) \wedge \big(\forall{b:B}\,\exists{a:A}\ f(a)=_V g(b)\big)\,,
    \]
    we put in a path in $V$ from $\mathsf{set}(A,f)$ to $\mathsf{set}(B,g)$.
    
  \item The ``set-truncation'' constructor: for all $x,y:V$ and all $p,q:x=_V y$, we add a (higher) path from $p$ to $q$.
\end{enumerate}
In (\ref{bisim}) we used the ``logical'' notation $\exists$ and $\forall$, etc., to indicate that we are working with the ``propositional truncations'' of the corresponding type theoretical operations $\Sigma$ and $\Pi$ (cf.~\cite{hottbook} ch.~3).  

Next, the membership relation $x\in y$ is defined for elements of $V$ by
\[
  (x \in \mathsf{set}(A,f))\ =_\mathrm{def}\  (\exists{a : A}.\ x =_V f(a))\,.
\]
One can then show entirely within the system that the resulting structure $(V,\in)$ satisfies \emph{almost} all of the axioms of Aczel's constructive set theory CZF \cite{AczelCZF} (e.g.\ Strong Collection is missing).

Finally, assuming the usual axiom of choice just for those types that are \emph{sets} in the sense of the homotopy levels, it then follows that $(V,\in)$ is a model of the full system of ZFC set theory (cf.~\cite{hottbook}, ch.~10.5).  The proofs of these results make essential use of the univalence axiom and have been fully formalized in the Coq proof assistant (cf.~\cite{formalizedVsets}).  

The system just mentioned is an interesting hybrid of classical set theory and constructive type theory that not only contains a model of ZFC but also many types of higher homotopy level that do not behave like classical sets.  This provides a new, more refined view of one possible relationship between classical and constructive foundations. Whereas constructive foundations are usually regarded as incompatible with classical logic, in this system the classical sets form a subsystem of constructive type theory consisting of certain objects distinguished by a natural, intrinsic, structural property---namely that their identity relation is always a proposition.

%%%%%%%%%%%%%%%%%%%%%%%%%%%%%%%%%%%%%%%%%%%%%%%%%%%%%%%%%%%%%%
%\subsection*{References and Further Information}
%
%More Information:
%\begin{center}
%\tt{www.HomotopyTypeTheory.org}\\
%\end{center}
%\bigskip
%
%The Book:
%\begin{center}
%\emph{Homotopy Type Theory: Univalent Foundations of Mathematics}\\
%The Univalent Foundations Program,\\
% Institute for Advanced Study, 
% Princeton, 2013
%\end{center}
%
%

\bibliographystyle{plain}

\bibliography{taitbib}

\end{document}